\documentclass[12pt]{amsart}%
\usepackage{graphicx}
\usepackage{amscd}
\usepackage{amsmath}%
\usepackage{amsfonts}%
\usepackage{amssymb}
\newcommand{\strutstretchdef}{\newcommand{\strutstretch}{\vphantom{\raisebox{1pt}{$\big($}\raisebox{-1pt}{$\big($}}}}
\theoremstyle{plain}
\newtheorem{theorem}{Theorem}[section]
\newtheorem{lemma}[theorem]{Lemma}

\newtheorem{corollary}[theorem]{Corollary}

\theoremstyle{definition}

\theoremstyle{remark}

\numberwithin{equation}{section}

\newlength{\struh}
\newlength{\textminustop}
\strutstretchdef
\begin{document}
\title[Weyl's Theorems and Local Spectral Theory]{\textbf{Weyl's Theorem, }$a$\textbf{-Weyl's theorem, }\\\textbf{and Local Spectral Theory} }
\author{Ra\'{u}l E. Curto}
\address{Department of Mathematics, The University of Iowa, Iowa City, Iowa 52242}
\email{curto@math.uiowa.edu}
\author{Young Min Han}
\address{Department of Mathematics, The University of Iowa, Iowa City, Iowa 52242}
\email{yhan@math.uiowa.edu}
\thanks{The research of the first named author was partially supported by NSF grants
DMS-9800931 and DMS-0099357. }

\begin{abstract}
We give necessary and sufficient conditions for a Banach space operator with
the single valued extension property (SVEP) to satisfy Weyl's theorem and
$a$-Weyl's theorem. We show that if $T$ or $T^{\ast}$ has SVEP and $T$ is
transaloid, then Weyl's theorem holds for $f(T)$ for every $f\in H(\sigma
(T))$. When $T^{\ast}$ has SVEP, $T$ is transaloid and $T$ is $a$-isoloid,
then $a$-Weyl's theorem holds for $f(T)$ for every $f\in H(\sigma(T))$. We
also prove that if $T$ or $T^{\ast}$ has SVEP, then the spectral mapping
theorem holds for the Weyl spectrum and for the essential approximate point spectrum.

\end{abstract}
\maketitle

\section{Introduction}

Let $B(X)$ denote the algebra of bounded linear operators acting on an
infinite dimensional complex Banach space $X$. If $T\in B(X)$ we shall write
$N(T)$ and $R(T)$ for the null space and range of $T$, respectively. \ Also,
let $\alpha(T):=\text{dim }N(T)$, $\beta(T):=\text{dim }X/R(T)$, and let
$\sigma(T)$ denote the spectrum of $T$, $\pi_{0}(T)$ the eigenvalues of $T$,
and $\pi_{0f}(T)$ the eigenvalues of finite multiplicity of $T$. \ An operator
$T\in B(X)$ is called \textit{Fredholm} if it has closed range, finite
dimensional null space, and its range has finite co-dimension. \ The
\textit{index} of a Fredholm operator $T$ is given by
\[
i(T):=\alpha(T)-\beta(T).
\]
An operator $T\in B(X)$ is called \textit{Weyl} if it is Fredholm of index
zero, and \textit{Browder} if it is Fredholm ``of finite ascent and descent;''
equivalently (\cite[Theorem 7.9.3]{Harte2}) if $T$ is Fredholm and $T-\lambda$
is invertible for sufficiently small $\lambda\neq0$ in $\mathbb{C}$. \ The
\textit{essential spectrum} $\sigma_{e}(T)$, the \textit{Weyl spectrum}
$\omega(T)$, and the \textit{Browder spectrum} $\sigma_{b}(T)$ of $T\in B(X)$
are defined by (\cite{Harte1}, \cite{Harte2})
\[
\sigma_{e}(T):=\{\lambda\in\mathbb{C}:T-\lambda\ \text{is not Fredholm}\},
\]%
\[
\omega(T):=\{\lambda\in\mathbb{C}:T-\lambda\ \text{is not Weyl}\},
\]
and
\[
\sigma_{b}(T):=\{\lambda\in\mathbb{C}:T-\lambda\ \text{is not Browder}\},
\]
respectively. \ Evidently
\[
\sigma_{e}(T)\subseteq\omega(T)\subseteq\sigma_{b}(T)=\sigma_{e}%
(T)\cup\text{acc}\,\sigma(T),
\]
where we write $\text{acc}\ K$ for the accumulation points of $K\subseteq
\mathbb{C}$.

If we write $\text{iso}\ K:=K\setminus\text{acc}\ K$ then we let
\[
\pi_{00}(T):=\{\lambda\in\text{iso}\ \sigma(T):0<\alpha(T-\lambda)<\infty\}
\]
denote the set of isolated eigenvalues of finite multiplicity, and we let
\[
p_{00}(T):=\sigma(T)\setminus\sigma_{b}(T)
\]
denote the set of Riesz points of $T$.

We say that Weyl's theorem holds for $T\in B(X)$\ if
\begin{equation}
\sigma(T)\setminus\omega(T)=\pi_{00}(T), \label{eq11}%
\end{equation}
and that Browder's theorem holds for $T\in B(X)$\ if
\begin{equation}
\sigma(T)\setminus\omega(T)=p_{00}(T). \label{eq12}%
\end{equation}
In \cite{Weyl}, H. Weyl proved that (\ref{eq11}) holds for hermitian
operators. Weyl's theorem has been extended from hermitian operators to
hyponormal operators, to Toeplitz operators \cite{Coburn}, and to several
classes of operators including seminormal operators (\cite{Berberian1},
\cite{Berberian2}).

In this article we give necessary and sufficient conditions for a Banach space
operator with the single valued extension property (SVEP) to satisfy Weyl's
theorem (Theorem \ref{thm22}) and $a$-Weyl's theorem (Corollary \ref{cor33}).
\ (For the relevant definitions please see below.) \ We show that if $T$ or
$T^{\ast}$ has SVEP and $T$ is transaloid, then Weyl's theorem holds for
$f(T)$ for every $f\in H(\sigma(T))$ (Theorems \ref{thm25} and \ref{thm35}).
We establish that if $T^{\ast}$ has SVEP, and if $T$ is transaloid and
$a$-isoloid, then $a$-Weyl's theorem holds for $f(T)$ for every $f\in
H(\sigma(T))$ (Theorem \ref{thm35}). We also prove that if $T$ or $T^{\ast}$
has SVEP, then the spectral mapping theorem holds for the Weyl spectrum
(Corollary \ref{cor26}) and for the essential approximate point spectrum
(Theorem \ref{thm31}).

To describe some of the above mentioned results, we consider the sets
\[
\Phi_{+}(X)=\{T\in B(X):R(T)\text{ is closed and}\ \alpha(T)<\infty\},
\]%
\[
\Phi_{-}(X)=\{T\in B(X):R(T)\text{ is closed and}\ \beta(T)<\infty\},
\]
and
\[
\Phi_{+}^{-}(X)=\{T\in B(X):T\in\Phi_{+}(X)\ \text{and}\ i(T)\leq0\}.
\]
Moreover,
\[
\sigma_{le}(T):=\{\lambda\in\mathbb{C}:T-\lambda\notin\Phi_{+}(X)\}
\]
is the \textit{left essential spectrum},
\[
\sigma_{re}(T):=\{\lambda\in\mathbb{C}:T-\lambda\notin\Phi_{-}(X)\}
\]
is the \textit{right essential spectrum},
\[
\sigma_{ea}(T):=\{\lambda\in\mathbb{C}:T-\lambda\notin\Phi_{+}^{-}(X)\}
\]
is the \textit{essential approximate point spectrum}, and
\[
\sigma_{a}(T):=\{\lambda\in\mathbb{C}:R(T-\lambda)\text{ is not closed or
}\alpha(T-\lambda)>0\}
\]
is the \textit{approximate point spectrum}. We now let
\[
\pi_{00}^{a}(T):=\{\lambda\in\text{iso}\ \sigma_{a}(T):0<\alpha(T-\lambda
)<\infty\ \},
\]
denote the set of eigenvalues of finite multiplicity which are isolated in
$\sigma_{a}(T)$, and we let%
\[
\sigma_{ab}(T):=\cap\{\sigma_{a}(T+K):TK=KT\ \text{and}\ K\in K(X)\}
\]
denote the \textit{Browder essential approximate point spectrum}, where $K(X)$
is the set of all compact operators on $X$. \ Observe that $\sigma
_{ea}(T)\subseteq\sigma_{ab}(T)$ for every $T\in B(X)$.

We say that $a$-Weyl's theorem holds for $T\in B(X)$ if\
\begin{equation}
\sigma_{a}(T)\setminus\sigma_{ea}(T)=\pi_{00}^{a}(T), \label{eq13}%
\end{equation}
and that $a$-Browder's theorem holds for $T\in B(X)$ if
\begin{equation}
\sigma_{ea}(T)=\sigma_{ab}(T). \label{eq14}%
\end{equation}
It is well known (\cite{Han}, \cite{Harte3}, \cite{Rakocevic1}) that if $T\in
B(X)$ then we have:
\[
\text{$a$-Weyl's theorem}\Longrightarrow\text{Weyl's theorem}\Longrightarrow
\text{Browder's theorem};
\]%
\[
\text{$a$-Weyl's theorem}\Longrightarrow\text{$a$-Browder's theorem}%
\Longrightarrow\text{Browder's theorem}.
\]

V. Rako\v{c}evi\'{c} \cite{Rakocevic1} has shown that (\ref{eq13}) holds for
cohyponormal operators. More recently, S.V. Djordjevi\'{c} and D.S.
Djordjevi\'{c} \cite{DD} have shown that if $T^{\ast}$ is quasihyponormal then
$a$-Weyl's theorem holds for $T$.

An operator $T\in B(X)$ is called \textit{isoloid} if every isolated point of
$\sigma(T)$ is an eigenvalue of $T$. If $T\in B(X)$, we write $r(T)$ for the
spectral radius of $T$; it is well known that $r(T)\leq||T||$. An operator
$T\in B(X)$ is called \textit{normaloid} if $r(T)=||T||$. $X\in B(X)$ is
called a \textit{quasiaffinity} if it has trivial kernel and dense range.
$S\in B(X)$ is said to be a \textit{quasiaffine transform} of $T\in B(X)$
(notation: $S\prec T$) if there is a quasiaffinity $X\in B(X)$ such that
$XS=TX$. If both $S\prec T$ and $T\prec S$, then we say that $S$ and $T$ are
\textit{quasisimilar}.

We say that $T\in B(X)$ has the \textit{single valued extension property}
(SVEP) if for every open set $U$ of $\mathbb{C}$ the only analytic solution
$f:U\longrightarrow X$ of the equation
\[
(T-\lambda)f(\lambda)=0\;\;\;(\lambda\in U)
\]
is the zero function (\cite{Foias}, \cite{Laursen3}). Given an arbitrary
operator $T\in B(X)$, the \textit{local resolvent set} $\rho_{T}(x)$ of $T$ at
the point $x\in X$ is defined as the union of all open subsets $U$ of
$\mathbb{C}$ for which there is an analytic function $f:U\longrightarrow X$
which satisfies
\[
(T-\lambda)f(\lambda)=x\ \quad\;(\text{all}\ \lambda\in U).
\]
The \textit{local spectrum} $\sigma_{T}(x)$ of $T$ at $x$ is then defined as
\[
\sigma_{T}(x):=\mathbb{C\setminus\rho}_{T}(x).
\]
For an arbitrary operator $T\in B(X)$, we define the \textit{local}\ (resp.
\textit{glocal}) \textit{spectral subspaces} of $T$ as follows. \ Given a set
$F\subseteq\mathbb{C}$ (resp. a closed set $G\subseteq\mathbb{C}$),
\[
X_{T}(F):=\{x\in X:\sigma_{T}(x)\subseteq F\}
\]
(resp.
\begin{align*}
\mathcal{X}_{T}(G)  &  :=\{x\in X:\text{there exists an analytic function }\\
f  &  :\mathbb{C\setminus}G\rightarrow X\text{ that satisfies }(T-\lambda
)f(\lambda)=x\text{ for all }\lambda\in\mathbb{C}\setminus G\}).
\end{align*}
An operator $T\in B(X)$ has \textit{Dunford's property} (C) if the local
spectral subspace $X_{T}(F)$ is closed for every closed set $F\subseteq
\mathbb{C}$; such operators automatically have SVEP \cite{Laursen3}.

We will use often the following two results.

\begin{lemma}
\label{lem11}(\cite[Theorem 2.6]{AM}) Let $T\in B(X)$ be an operator
satisfying SVEP, and let $\lambda\notin\sigma_{e}(T)$. \ Then $T-\lambda$ has
SVEP if and only if $T-\lambda$ has finite ascent.
\end{lemma}

\begin{lemma}
\label{lem12}(\cite[Theorem 3.1]{Koliha}) Let $T\in B(X)$ and assume that
$\lambda\in\pi_{00}(T)$. Then $T=T_{1}\oplus T_{2}$ with respect to the
decomposition $X=\mathcal{X}_{T}(\{\lambda\})\oplus K(T-\lambda)$, and
$\sigma(T_{1})=\{\lambda\}$, $\sigma(T_{2})=\sigma(T)\setminus\{\lambda\}$,
where $K(T):=\{x\in X:Tx_{n+1}=x_{n},Tx_{1}=x,$ $||x_{n}||\leq c^{n}%
||x||\ (n=1,2,\ldots)\ \text{for some}\ c>0,x_{n}\in X\}$.
\end{lemma}

\section{Extensions of Weyl's theorem}

The following theorem relates Weyl's theorem to local spectral theory. \ As
motivation for the proof, we use some ideas in \cite{Laursen2} and
\cite{Laursen3}.

\begin{theorem}
Let $T\in B(X)$ and assume that $\mathcal{X}_{T}(\{\lambda\})$ is finite
dimensional for each $\lambda\in\pi_{0f}(T)$. \ Then Weyl's theorem holds for
$T$.
\end{theorem}

\begin{proof}
We must show that $\sigma(T)\setminus\omega(T)=\pi_{00}(T)$. \ Suppose that
$\lambda\in\sigma(T)\setminus\omega(T)$. \ Then $T-\lambda$ is Weyl but not
invertible. \ Then $\lambda\in\pi_{0f}(T)$, and hence dim $\mathcal{X}%
_{T}(\{\lambda\})<\infty$. \ Since $\mathcal{X}_{T}(\{\lambda\})$ is an
invariant subspace for $T$, we can write
\[
T=\left(
\begin{array}
[c]{cc}%
A & C\\
0 & B
\end{array}
\right)  ,
\]
where $A:=T|_{\mathcal{X}_{T}(\{\lambda\})}$. \ Since dim $\mathcal{X}%
_{T}(\{\lambda\})<\infty$, $T-\lambda$ is Weyl if and only if $B-\lambda$ is
Weyl. \ Thus, $B-\lambda$ is Weyl. \ However, $N(T-\lambda)\subseteq
\mathcal{X}_{T}(\{\lambda\})$, hence $B-\lambda$ is injective. \ Therefore
$B-\lambda$ is invertible. \ Recall now that since dim $\mathcal{X}%
_{T}(\{\lambda\})<\infty$, $\sigma(A)$ must be finite and $\sigma(T)=\sigma(A)%
{\textstyle\bigcup}
\sigma(B)$. \ Thus, $\lambda$ is an isolated point of $\sigma(T)$%
.\ \ Therefore $\lambda\in$ iso $\sigma(T)\setminus\omega(T)$, and by the
punctured neighborhood theorem, $\lambda\in\pi_{00}(T)$. \ \newline
Conversely, suppose that $\lambda\in\pi_{00}(T)$. By Lemma \ref{lem12},
$T=T_{1}\oplus T_{2}$ on $\mathcal{X}_{T}(\{\lambda\})\oplus K(T-\lambda)$,
with $\sigma(T_{1})=\{\lambda\}$ and $\sigma(T_{2})=\sigma(T)\setminus
\{\lambda\}$. \ Since dim $\mathcal{X}_{T}(\{\lambda\})<\infty$, $T-\lambda$
is Weyl. Therefore $\lambda\in\sigma(T)\setminus\omega(T)$.
\end{proof}

In \cite{Jeon}, it was shown that if $T\in B(X)$ has \textit{totally finite
ascent} (in the sense that $T-\lambda$ has finite ascent for each $\lambda
\in\mathbb{C}$), then Weyl's theorem holds for $T$ if and only if
$R(T-\lambda)$ is closed for all $\lambda\in\pi_{00}(T)$. \ In general, if $T$
has totally finite ascent then $T$ has SVEP(\cite{Laursen1}). However, the
converse is not true. \ Consider the following example: let $T\in B(l_{2})$ be
given by
\begin{equation}
T(x_{0},x_{1},x_{2},\cdots):=({\frac{1}{2}}x_{1},{\frac{1}{3}}x_{2},\cdots).
\label{operator}%
\end{equation}
Then clearly $T$ does not have finite ascent. \ But since $T$ is
quasinilpotent, $T$ has SVEP. \ Thus SVEP is a much weaker condition than
having totally finite ascent. \ However, we can prove:

\begin{theorem}
\label{thm22}Suppose that $T\in B(X)$ has SVEP. \ Then the following
statements are equivalent:\newline (i) Weyl's theorem holds for $T$;
\newline (ii) $R(T-\lambda)$ is closed for all $\lambda\in\pi_{00}(T)$;
\newline (iii) $\mathcal{X}_{T}(\{\lambda\})$ is finite dimensional for every
$\lambda\in\pi_{00}(T)$; \newline (iv) $\gamma_{T}(\zeta)$ is discontinuous on
$\pi_{00}(T)$, where $\gamma_{T}(\cdot)$ denotes the reduced minimum modulus
of $T$, i.e.,
\[
\gamma_{T}(\zeta):=\inf\{\frac{||(T-\zeta)x||}{\text{dist}(x,N(T-\zeta))}:x\in
X\setminus N(T-\zeta)\}.
\]
\end{theorem}

\begin{proof}
(i)$\Rightarrow$(ii): Suppose $\lambda\in\pi_{00}(T)$. Since Weyl's theorem
holds for $T$, $\lambda\in\sigma(T)\setminus\omega(T)$. Therefore
$R(T-\lambda)$ is closed.\newline (ii)$\Rightarrow$(iii): Let $\lambda\in
\pi_{00}(T)$. By Lemma \ref{lem12}, $T=T_{1}\oplus T_{2}$ on $\mathcal{X}%
_{T}(\{\lambda\})\oplus K(T-\lambda)$, where $\sigma(T_{1})=\{\lambda\}$ and
$\sigma(T_{2})=\sigma(T)\setminus\{\lambda\}$. Since $R(T-\lambda)$ is closed,
an application of the punctured neighborhood theorem shows that $\lambda$ is a
Riesz point of $T$. \ Therefore $R(P)$ is finite dimensional, where $P\in
B(X)$ is the spectral projection corresponding to $\lambda$, given by
\[
P:=\frac{1}{2\pi i}\ {\int_{\partial D}}(z-T)^{-1}dz,
\]
where $D$ is an open disk of center $\lambda$ which contains no other points
of $\sigma(T)$. \ Now,
\[
R(P)=\{x\in X:\lim_{n\rightarrow\infty}||(T-\lambda)^{n}x||^{\frac{1}{n}%
}=0\}=\mathcal{X}_{T}(\{\lambda\});
\]
hence $\mathcal{X}_{T}(\{\lambda\})$ is finite dimensional.\newline
(iii)$\Rightarrow$(i): Suppose $\lambda\in\sigma(T)\setminus\omega(T)$. \ Then
$T-\lambda$ is Weyl but not invertible. \ We first show that $\lambda
\in\partial\sigma(T)$. \ Assume to the contrary that $\lambda\in\text{int
}\sigma(T)$. \ Then there exists a neighborhood $U$ of $\lambda$ such that
$\text{dim}\ N(T-\mu)>0$ for all $\mu\in U$. \ It follows from \cite[Theorem
10]{Finch} that $T$ does not have SVEP, a contradiction. \ Therefore
$\lambda\in\partial\sigma(T)\setminus\omega(T)$, and it follows from the
punctured neighborhood theorem that $\lambda\in\pi_{00}(T)$. \ Conversely,
suppose that $\lambda\in\pi_{00}(T)$. \ By Lemma \ref{lem12}, $T=T_{1}\oplus
T_{2}$ on $\mathcal{X}_{T}(\{\lambda\})\oplus K(T-\lambda)$, where
$\sigma(T_{1})=\{\lambda\}$ and $\sigma(T_{2})=\sigma(T)\setminus\{\lambda\}$.
\ Since $\mathcal{X}_{T}(\{\lambda\})$ is a finite dimensional subspace of
$X$, $T-\lambda$ is Weyl. \ Therefore $\lambda\in\sigma(T)\setminus\omega
(T)$.\newline (i)$\Leftrightarrow$(iv): If Weyl's theorem holds for $T$, then
it follows from \cite[Theorem 1]{Gustafson} that $\gamma_{T}(\lambda)$ is
discontinuous for each $\lambda\in\pi_{00}(T)$. \ Conversely, suppose
$\gamma_{T}(\lambda)$ is discontinuous on $\pi_{00}(T)$. \ Since $T$ has SVEP,
it suffices to show that $\pi_{00}(T)\subseteq\sigma(T)\setminus\omega(T)$.
\ Suppose that $\lambda\in\pi_{00}(T)$. \ Since $\lambda$ is an isolated point
of $\sigma(T)$, there exist $\delta>0$ and an open disk $D(\lambda,\delta)$
centered in $\lambda$ such that $D(\lambda,\delta)\cap\sigma(T)=\{\lambda\}$.
\ For every $\mu\in D(\lambda,\delta)\setminus\{\lambda\}$, we have $T-\mu$
injective and thus,
\begin{align*}
\gamma_{T}(\mu) &  =\inf_{x\in X\setminus N(T-\mu)}\frac{||(T-\mu
)x||}{\text{dist}(x,N(T-\mu))}=\inf_{x\neq0}\frac{||(T-\mu)x||}{||x||}\\
&  \leq\inf_{x\in N(T-\lambda)\setminus\{0\}}\frac{||(T-\lambda)x-(\mu
-\lambda)x||}{||x||}\\
&  =\inf_{x\in N(T-\lambda)\setminus\{0\}}\frac{||(\mu-\lambda)x||}%
{||x||}=|\mu-\lambda|.
\end{align*}
Since $\gamma_{T}(\zeta)$ is discontinuous at $\lambda$, $\gamma_{T}%
(\lambda)>0$. \ Therefore $R(T-\lambda)$ is closed, and hence $T-\lambda$ is
Weyl. \ This completes the proof.
\end{proof}

Before we state our next theorem, we need a definition and two preliminary
results. \ Recall that an operator $T\in B(X)$ is called \textit{transaloid}%
\ if $T-\lambda$ is normaloid for every $\lambda\in\mathbb{C}$.

\begin{lemma}
\label{lem23}Suppose that $T\in B(X)$ is transaloid. \ Then
\[
\mathcal{X}_{T}(\{\lambda\})=N(T-\lambda)\quad\text{for every}\ \lambda
\in\mathbb{C}.
\]
\end{lemma}

\begin{proof}
Observe that $\mathcal{X}_{T}(\{\lambda\})=\{x\in X:\lim\limits_{n\rightarrow
\infty}||(T-\lambda)^{n}x||^{\frac{1}{n}}=0\}$ for each $\lambda\in\mathbb{C}%
$. \ Since $T$ is transaloid, $T-\mu$ is normaloid for each $\mu\in\mathbb{C}%
$. \ Therefore $||(T-\lambda)x||\leq||(T-\lambda)^{n}x||^{\frac{1}{n}}$ for
all $x\in X$ and $n\in\mathbb{N}$, and hence $\mathcal{X}_{T}(\{\lambda
\})\subseteq N(T-\lambda)$ for every $\lambda\in\mathbb{C}$. \ The converse is clear.
\end{proof}

\begin{lemma}
\label{lem24}Suppose that $T\in B(X)$ has SVEP and is transaloid. \ Then $T$
is isoloid.
\end{lemma}

\begin{proof}
Suppose that $\lambda$ is an isolated point of $\sigma(T)$. \ Then it follows
from \cite{Laursen3} that $X=\mathcal{X}_{T}(\{\lambda\})+X_{T}%
(\mathbb{C\setminus\{\lambda\})}$. \ Assume to the contrary that $T-\lambda$
is injective. \ It follows from Lemma \ref{lem23} that $X=X_{T}%
(\mathbb{C\setminus\{\lambda\})}$. \ But $(T-\lambda)X_{T}(\mathbb{C\setminus
\{\lambda\})}=X_{T}\mathbb{(C\setminus\{\lambda\})}$; hence $T-\lambda$ is
surjective. \ Therefore $T-\lambda$ is invertible, a contradiction. \ It
follows that $T$ is isoloid.
\end{proof}

In the following theorem, recall that $H(\sigma(T))$ is the space of functions
analytic in an open neighborhood of $\sigma(T)$.

\begin{theorem}
\label{thm25}Suppose that $T\in B(X)$ has SVEP and is transaloid. \ Then
Weyl's theorem holds for $f(T)$ for every $f\in H(\sigma(T))$.
\end{theorem}

\begin{proof}
We first show that Weyl's theorem holds for $T$. \ Since SVEP and being
transaloid are translation-invariant properties, it suffices to show that
\[
0\in\pi_{00}(T)\Longleftrightarrow T\ \text{is Weyl and not invertible}.
\]
Suppose $0\in\pi_{00}(T)$. \ Then using the spectral projection $P:=\frac{1}%
{2\pi i}\int_{\partial D}(\lambda-T)^{-1}d\lambda$, where $D$ is an open disk
of center $0$ which contains no other points of $\sigma(T)$, we can write
$T=T_{1}\oplus T_{2}$, where $\sigma(T_{1})=\{0\}$ and $\sigma(T_{2}%
)=\sigma(T)\setminus\{0\}$. \ It follows from Lemma \ref{lem23} that
$P(X)=\{x\in X:\lim\limits_{n\rightarrow\infty}||T^{n}x||^{\frac{1}{n}%
}=0\}=\mathcal{X}_{T}(\{0\})=N(T)$. \ Since $N(T)$ is a finite dimensional
subspace of $X$, we must have $\omega(T)=\omega(T_{2})$. \ But $T_{2}$ is
invertible, hence $T$ is Weyl. Therefore $0\in\sigma(T)\setminus\omega(T)$.
\ Conversely, suppose that $0\in\sigma(T)\setminus\omega(T)$. \ Then it
follows from Lemma \ref{lem23} that $\mathcal{X}_{T}(\{0\})=N(T)$. \ Since
$\mathcal{X}_{T}(\{0\})$ is a closed invariant subspace for $T$, $T$ can be
represented by the following $2\times2$ operator matrix:
\[
T=%
\begin{pmatrix}
0 & T_{1}\\
0 & T_{2}%
\end{pmatrix}
.
\]
Since $\mathcal{X}_{T}(\{0\})$ is a finite dimensional subspace of $X$, $T$ is
Weyl if and only if $T_{2}$ is Weyl. \ But $\mathcal{X}_{T}(\{0\})=N(T)$;
hence $T_{2}$ is injective, and so $T_{2}$ is invertible. \ It follows from
the punctured neighborhood theorem that $0\in\pi_{00}(T)$. \ Thus Weyl's
theorem holds for $T$. \ \newline We now claim that
\begin{equation}
f(\omega(T))=\omega(f(T))\text{ for all }f\in H(\sigma(T)).\label{SMTWeyl}%
\end{equation}
\ Let $f\in H(\sigma(T))$. \ Since $\omega(f(T))\subseteq f(\omega(T))$ with
no restriction on $T$, it suffices to show that $f(\omega(T))\subseteq
\omega(f(T))$. \ Suppose $\lambda\notin\omega(f(T))$. \ Then $f(T)-\lambda$ is
Weyl and
\begin{equation}
f(T)-\lambda=c(T-\alpha_{1})(T-\alpha_{2})\cdots(T-\alpha_{n})g(T),\label{251}%
\end{equation}
where $c,\alpha_{1},\alpha_{2},\cdots,\alpha_{n}\in\mathbb{C}$, and $g(T)$ is
invertible. \ Since the operators on the right-hand side of (\ref{251})
commute, every $T-\alpha_{i}$ is Fredholm. \ Since $T$ has SVEP, it follows
from Lemma \ref{lem11} that each $T-\alpha_{i}$ has finite ascent. \ Now we
show that $i(T-\alpha_{i})\leq0$ for each $i=1,2,\cdots,n$. \ Observe that if
$A\in B(X)$ is Fredholm of finite ascent then $i(A)\leq0$ : indeed, either $A$
has finite descent, in which case $A$ is Browder and $i(A)=0$, or $A$ has
infinite descent and
\[
n\cdot i(A)=\alpha(A^{n})-\beta(A^{n})\longrightarrow-\infty\quad
\text{as}\ n\longrightarrow\infty,
\]
which implies that $i(A)<0$. \ Thus $i(T-\alpha_{i})\leq0$ for each
$i=1,2,\cdots,n$. \ Therefore $\lambda\notin f(\omega(T))$, and hence
$f(\omega(T))=\omega(f(T))$. \ \newline We now recall that if $T$ is isoloid
then
\begin{equation}
f(\sigma(T)\setminus\pi_{00}(T))=\sigma(f(T))\setminus\pi_{00}%
(f(T))\label{isoloid}%
\end{equation}
$\text{for every}\ f\in H(\sigma(T))\ $\cite[Lemma]{Lee}. By Lemma
\ref{lem24}, $T$ is isoloid, and since Weyl's theorem holds for $T$, we have%
\begin{align*}
\sigma(f(T))\setminus\pi_{00}(f(T)) &  =f(\sigma(T)\setminus\pi_{00}%
(T))\;\;\text{(by(\ref{isoloid})}\\
&  =f(\omega(T))=\omega(f(T))\;\;\text{(by \ref{SMTWeyl}).}%
\end{align*}
Therefore
\[
\pi_{00}(f(T))=\sigma(f(T))\setminus\lbrack\sigma(f(T))\setminus\pi
_{00}(f(T))]=\sigma(f(T))\setminus\omega(f(T)),
\]
so Weyl's theorem holds for $f(T)$, as desired.
\end{proof}

From the proof of Theorem \ref{thm25} we obtain the following useful consequence.

\begin{corollary}
\label{cor26}Let $T\in B(X)$. \ Suppose that $T$ or $T^{\ast}$ has SVEP.
\ Then
\[
\omega(f(T))=f(\omega(T))\quad\text{for every}\ f\in H(\sigma(T)).
\]
\end{corollary}

\section{Extensions of $a$-Weyl's theorem}

Let $T\in B(X)$. It is known that the inclusion $\sigma_{ea}(f(T))\subseteq
f(\sigma_{ea}(T))$ holds for every $f\in H(\sigma(T))$, with no restriction on
$T$ \cite{Rakocevic2}. The next theorem shows that the spectral mapping
theorem holds for the essential approximate point spectrum, for operators
having SVEP.

\begin{theorem}
\label{thm31}Let $T\in B(X)$ and suppose that $T$ or $T^{\ast}$ has SVEP.
\ Then
\[
\sigma_{ea}(f(T))=f(\sigma_{ea}(T))\quad\text{for every}\ f\in H(\sigma(T)).
\]
\end{theorem}

\begin{proof}
Let $f\in H(\sigma(T))$. \ It suffices to show that $f(\sigma_{ea}%
(T))\subseteq\sigma_{ea}(f(T))$. \ Suppose that $\lambda\notin\sigma
_{ea}(f(T))$. \ Then $f(T)-\lambda\in\Phi_{+}^{-}(X)$ and
\begin{equation}
f(T)-\lambda=c(T-\alpha_{1})(T-\alpha_{2})\cdots(T-\alpha_{n})g(T),
\label{311}%
\end{equation}
where $c,\alpha_{1},\alpha_{2},\cdots,\alpha_{n}\in\mathbb{C}$, and $g(T)$ is
invertible. \ Since the operators on the right-hand side of (\ref{311})
commute, $T-\alpha_{i}\in\Phi_{+}(X)$. \ Since $T$ has SVEP, it follows from
Lemma \ref{lem11} that each $T-\alpha_{i}$ has finite ascent. \ Therefore by
the proof of Theorem \ref{thm25}, $i(T-\alpha_{i})\leq0$ for each
$i=1,2,\cdots,n$. \ It follows that $\lambda\notin f(\sigma_{ea}(T))$.
\ \newline Suppose now that $T^{\ast}$ has SVEP. \ Since $T-\alpha_{i}\in
\Phi_{+}(X)$, $T^{\ast}-\alpha_{i}\in\Phi_{-}(X^{\ast})$. \ Since $T^{\ast}$
has SVEP, it follows from Lemma \ref{lem11} that each $T-\alpha_{i}$ has
finite descent. \ We claim that $i(T-\alpha_{i})\geq0$ for each $i=1,2,\cdots
,n$. \ Observe that if $A\in\Phi_{-}(X)$ and $A$ is not Fredholm then
evidently $i(A)\geq0$. \ If $A$ is Fredholm with finite descent, then either
$A$ has finite ascent (and then $A$ is Browder and $i(A)=0$), or $A$ has
infinite ascent (and then
\[
n\cdot i(A)=\alpha(A^{n})-\beta(A^{n})\longrightarrow\infty\quad
\text{as}\ n\longrightarrow\infty,
\]
which implies that $i(A)>0$). \ Thus $i(T-\alpha_{i})\geq0$ for each
$i=1,2,\cdots,n$. \ However,
\[
0\leq\sum_{i=1}^{n}i(T-\alpha_{i})=i(f(T)-\lambda)\leq0,
\]
and so $T-\alpha_{i}$ is Weyl for each $i=1,2,\cdots,n$. \ Hence
$\lambda\notin f(\sigma_{ea}(T))$, and so $\sigma_{ea}(f(T))=f(\sigma
_{ea}(T))$. \ This completes the proof of the theorem.
\end{proof}

In general, we cannot expect that $a$-Weyl's theorem necessarily holds for
operators having SVEP. \ Consider the quasinilpotent operator $T$ on $\ell
_{2}$ given by (\ref{operator}); \ this operator has SVEP. \ But $\sigma
_{a}(T)=\sigma_{ea}(T)=\{0\}$, and $\pi_{00}^{a}(T)=\{0\}$; hence $a$-Weyl's
theorem does not hold for $T$. \ However, $a$-Browder's theorem does hold, as
the following result shows.

\begin{theorem}
\label{thm32}Suppose that $T\in B(X)$ has SVEP and $S\in B(X)$ satisfies
$S\prec T$. \ Then $a$-Browder's theorem holds for $f(S)$, for every $f\in
H(\sigma(S))$.
\end{theorem}

\begin{proof}
We first recall that $S$ has SVEP. \ For, let $U$ be any open set and
$f:U\longrightarrow X$ be any analytic function such that $(S-\lambda
)f(\lambda)=0$ for all $\lambda\in U$. \ Since $S\prec T$, there exists a
quasiaffinity $A$ such that $AS=TA$. \ Thus, $A(S-\lambda)=(T-\lambda)A$ for
all $\lambda\in U$. \ Since $(S-\lambda)f(\lambda)=0$ for all $\lambda\in U$,
$0=A(S-\lambda)f(\lambda)=(T-\lambda)Af(\lambda)$ for all $\lambda\in U$.
\ But $T$ has SVEP, hence $Af(\lambda)=0$ for all $\lambda\in U$. \ Since $A$
is a quasiaffinity, $f(\lambda)=0$ for all $\lambda\in U$. \ Therefore $S$ has
SVEP. \ Next we show that $a$-Browder's theorem holds for $S$, i.e., that
$\sigma_{ea}(S)=\sigma_{ab}(S)$. \ It is well known that $\sigma
_{ea}(S)\subseteq\sigma_{ab}(S)$. \ To prove the converse, suppose that
$\lambda\in\sigma_{a}(S)\setminus\sigma_{ea}(S)$. \ Then $S-\lambda\in\Phi
_{+}^{-}(X)$ and $S-\lambda$ is not bounded below. \ Since $S$ has SVEP and
$S-\lambda\in\Phi_{+}^{-}(X)$, it follows from Lemma \ref{lem11} that
$S-\lambda$ has finite ascent. Therefore by \cite[Theorem 2.1]{Rakocevic2},
$\lambda\in\sigma_{a}(S)\setminus\sigma_{ab}(S)$. \ Thus $a$-Browder's theorem
holds for $S$. \ Hence, it follows from Theorem \ref{thm31} that
\[
\sigma_{ab}(f(S))=f(\sigma_{ab}(S))=f(\sigma_{ea}(S))=\sigma_{ea}(f(S))
\]
(all $f\in H(\sigma(S))$), and so $a$-Browder's theorem holds for $f(S)$.
\end{proof}

In analogy with Theorem \ref{thm22}, we obtain

\begin{corollary}
\label{cor33}Suppose that $T\in B(X)$ has SVEP and $S\in B(X)$ satisfies
$S\prec T$. \ The following statements are equivalent:\newline (i) $a$-Weyl's
theorem holds for $S$;\newline (ii) $R(S-\lambda)$ is closed for all
$\lambda\in\pi_{00}^{a}(S)$;\newline (iii) $\gamma_{S}(\zeta)$ is
discontinuous on $\pi_{00}^{a}(S)$, where $\gamma_{S}(\cdot)$ denotes the
reduced minimum modulus;\newline (iv) $\sigma_{ea}(S)\cap\pi_{00}%
^{a}(S)=\emptyset$;\newline (v) $\pi_{00}^{a}(S)=\sigma_{a}(S)\setminus
\sigma_{ab}(S)$.
\end{corollary}

\begin{proof}
Since $T$ has SVEP and $S\prec T$, it follows from Theorem \ref{thm32} that
$a$-Browder's theorem holds for $S$. \ Therefore $\sigma_{ea}(S)=\sigma
_{ab}(S)$. \newline (i)$\Leftrightarrow$(ii): Suppose $\lambda\in\pi_{00}%
^{a}(S)$. Since $a$-Weyl's theorem holds for $S$, $\lambda\in\sigma
_{a}(S)\setminus\sigma_{ea}(S)$. Therefore $R(S-\lambda)$ is closed.
\ Conversely, suppose that $R(S-\lambda)$ is closed for all $\lambda\in
\pi_{00}^{a}(S)$. \ Since $a$-Browder's theorem holds for $S$, $\sigma
_{ea}(S)=\sigma_{ab}(S)$. \ It follows from \cite[Corollary 2.2]{Rakocevic2}
that $\sigma_{a}(S)\setminus\sigma_{ea}(S)=\sigma_{a}(S)\setminus\sigma
_{ab}(S)\subseteq\pi_{00}^{a}(S)$. \ Conversely, let $\lambda\in\pi_{00}%
^{a}(S)$. \ Then $S-\lambda$ has closed range, and so $S-\lambda\in\Phi
_{+}(X)$. \ Since $S$ has SVEP, it follows from Lemma \ref{lem11} that
$S-\lambda$ has finite ascent. \ Therefore $i(S-\lambda)\leq0$, and hence
$S-\lambda\in\Phi_{+}^{-}(X)$. \newline (i)$\Leftrightarrow$(iii): If
$a$-Weyl's theorem holds for $S$, then it follows from \cite[Theorem
2.4]{Rakocevic1} that $\gamma_{S}(\lambda)$ is discontinuous for each
$\lambda\in\pi_{00}^{a}(S)$. \ Conversely, suppose that $\gamma_{S}(\lambda)$
is discontinuous on $\pi_{00}^{a}(S)$. \ To show that $a$-Weyl's theorem holds
for $S$, it suffices to show that $R(S-\lambda)$ is closed for all $\lambda
\in\pi_{00}^{a}(S)$. \ Let $\lambda\in\pi_{00}^{a}(S)$. \ Since $\lambda$ is
an isolated point of $\sigma_{a}(S)$, there exist $\epsilon>0$ and an open
disk $D(\lambda,\epsilon)$ centered in $\lambda$ such that $D(\lambda
,\epsilon)\cap\sigma_{a}(S)=\{\lambda\}$. \ For every $\mu\in D(\lambda
,\epsilon)\setminus\{\lambda\}$, we have $S-\mu$ injective, and therefore
\begin{align*}
\gamma_{S}(\mu) &  =\inf_{x\in X\setminus N(S-\mu)}\frac{||(S-\mu
)x||}{\text{dist}(x,N(S-\mu))}=\inf_{x\neq0}\frac{||(S-\mu)x||}{||x||}\\
&  \leq\inf_{x\in N(S-\lambda)\setminus\{0\}}\frac{||(S-\lambda)x-(\mu
-\lambda)x||}{||x||}\\
&  =\inf_{x\in N(S-\lambda)\setminus\{0\}}\frac{||(\mu-\lambda)x||}%
{||x||}=|\mu-\lambda|.
\end{align*}
Since $\gamma_{S}(\zeta)$ is discontinuous at $\lambda$, $\gamma_{S}%
(\lambda)>0$. \ It follows that $R(S-\lambda)$ is closed. \newline
(ii)$\Leftrightarrow$(iv): Assume to the contrary that $\lambda\in\sigma
_{ea}(S)\cap\pi_{00}^{a}(S)$, so $R(S-\lambda)$ is closed. \ Since $S$ has
SVEP, $i(S-\lambda)\leq0$. Therefore $S-\lambda\in\Phi_{+}^{-}(X)$, and so
$\lambda\notin\sigma_{ea}(S)$, a contradiction. \ Conversely, let $\lambda
\in\pi_{00}^{a}(S)$. \ Since $\sigma_{ea}(S)\cap\pi_{00}^{a}(S)=\emptyset$,
$\lambda\notin\sigma_{ea}(S)$. \ Therefore $R(S-\lambda)$ is closed.
\newline (i)$\Leftrightarrow$(v): Since $a$-Weyl's theorem holds for $S$, we
have
\[
\pi_{00}^{a}(S)=\sigma_{a}(S)\setminus\sigma_{ea}(S)=\sigma_{a}(S)\setminus
\sigma_{ab}(S).
\]
Conversely, suppose that $\pi_{00}^{a}(S)=\sigma_{a}(S)\setminus\sigma
_{ab}(S)$. \ Since $S$ has SVEP, $\sigma_{ea}(S)=\sigma_{ab}(S)$. \ Therefore
$\pi_{00}^{a}(S)=\sigma_{a}(S)\setminus\sigma_{ea}(S)$, and hence $a$-Weyl's
theorem holds for $S$. \ This completes the proof.
\end{proof}

Recall \cite[Definition 13]{Harte3} that an operator $T\in B(X)$ is called
\textit{reguloid} if each isolated point of its spectrum is a regular point,
in the sense that there is a \textit{generalized inverse} $S_{\lambda}\in
B(X)$, i.e., $(T-\lambda)=(T-\lambda)S_{\lambda}(T-\lambda)$.

\begin{corollary}
\label{cor34}Let $T\in B(X)$. Suppose that $T$ or $T^{\ast}$ has Dunford's
property (C) and $T$ is reguloid. Then Weyl's theorem holds for $f(T)$ for
every $f\in H(\sigma(T))$.
\end{corollary}

\begin{proof}
Since Dunford's property (C) implies SVEP \cite[Proposition 1.2.19]{Laursen3},
it follows from Theorem \ref{thm32} that Browder's theorem holds for $T$.
\ But $T$ is reguloid, hence $T-\lambda$ has closed range for each $\lambda
\in\pi_{00}(T)$. \ Therefore Weyl's theorem holds for $T$ by Theorem
\ref{thm22}. \ Since $T$ is reguloid, it is also isoloid by \cite[Theorem
14]{Harte3}. \ Hence by the proof of Theorem \ref{thm25}, Weyl's theorem holds
for $f(T)$, for each $f\in H(\sigma(T))$.
\end{proof}

Finally, recall that $T\in B(X)$ is called \textit{approximate isoloid}%
\ ($a$-isoloid) if every isolated point of $\sigma_{a}(T)$ is an eigenvalue of
$T$. \ We observe that if $T$ is $a$-isoloid then it is isoloid because the
boundary of $\sigma(T)$ is contained in $\sigma_{a}(T)$. \ (However, the
converse is not true. \ Consider the following example: let $T=T_{1}\oplus
T_{2}$, where $T_{1}$ is the unilateral shift on\ $l_{2}$ and $T_{2}$ is
injective and quasinilpotent on $l_{2}$. \ Then $\sigma(T)=\{z\in
C:|z|\leq1\}$ \text{and} $\sigma_{a}(T)=\{z\in C:|z|=1\}\cup\{0\}$.
\ Therefore $T$ is isoloid but not $a$-isoloid.)

\begin{theorem}
\label{thm35}Suppose that $T^{\ast}\in B(X^{\ast})$ has SVEP and $T$ is
transaloid. \ Then $a$-Weyl's theorem holds for $T$. \ If, in addition, $T$ is
$a$-isoloid, then $a$-Weyl's theorem holds for $f(T)$, for every $f\in
H(\sigma(T))$.
\end{theorem}

\begin{proof}
We first show that $a$-Weyl's theorem holds for $T$. \ Since SVEP and being
transaloid are translation-invariant properties, it suffices to show that
\[
0\in\pi_{00}^{a}(T)\Longleftrightarrow0\in\sigma_{a}(T)\setminus\sigma
_{ea}(T).
\]
Suppose that $0\in\pi_{00}^{a}(T)$. \ Since $T^{\ast}$ has SVEP, it follows
from that [8, Corollary 7] that $\sigma(T)=\sigma_{a}(T)$. \ Therefore $0$ is
an isolated point of $\sigma(T)$. \ Using the spectral projection
$P:=\frac{1}{2\pi i}\int_{\partial D}(\lambda-T)^{-1}d\lambda$, where $D$ is
an open disk of center $0$ which contains no other points of $\sigma(T)$, we
can represent $T$ as the direct sum
\[
T=%
\begin{pmatrix}
T_{1} & 0\\
0 & T_{2}%
\end{pmatrix}
,
\]
where$\ \sigma(T_{1})=\{0\}\ $and$\ \sigma(T_{2})=\sigma(T)\setminus\{0\}$.
$\ $It follows from Lemma \ref{lem23} that
\[
P(X)=\{x\in X:\lim_{n\rightarrow\infty}||T^{n}x||^{\frac{1}{n}}%
=0\}=\mathcal{X}_{T}(\{0\})=N(T).
\]
Since $N(T)$ is a finite dimensional subspace of $X$, $\omega(T)=\omega
(T_{2})$. \ But $T_{2}$ is invertible, hence $T$ is Weyl. Therefore
$0\in\sigma_{a}(T)\setminus\sigma_{ea}(T)$. \ Conversely, suppose that
$0\in\sigma_{a}(T)\setminus\sigma_{ea}(T)$. \ Then $T\in\Phi_{+}(X)$ and
$i(T)\leq0$. \ Since $T^{\ast}$ has SVEP, it follows from Lemma 1.1 that
$i(T)\geq0$. Therefore $T$ is Weyl, and so $0\in\sigma_{a}(T)\setminus
\omega(T)$. \ Observe now that Browder's theorem holds for $T^{\ast}$ by
Theorem \ref{thm32}. Also, Browder's theorem holds for $T^{\ast}$ if and only
if it holds for $T$. \ Thus Browder's theorem holds for $T$, and so $0\in
\pi_{00}^{a}(T)$. \ \newline We have thus established that $a$-Weyl's theorem
holds for $T$. \ Let $f\in H(\sigma(T))$. \ Since $T^{\ast}$ has SVEP, it
follows from Theorem \ref{thm31} that $\sigma_{ea}(f(T))=f(\sigma_{ea}(T))$.
\ But $a$-Weyl's theorem holds for $T$, hence $\sigma_{ea}(T)=\sigma_{ab}(T)$.
\ Therefore
\[
\sigma_{ea}(f(T))=f(\sigma_{ea}(T))=f(\sigma_{ab}(T))=\sigma_{ab}(f(T))
\]
(the last equality by virtue of \cite[Theorem 3.4]{Rakocevic2}) and so
$a$-Browder's theorem holds for $f(T)$. \ Now let $\lambda\in\pi_{00}%
^{a}(f(T))$. \ Since $\sigma(T)=\sigma_{a}(T)$, it follows that $\sigma
(f(T)=\sigma_{a}(f(T))$, so $\lambda\in\pi_{00}(f(T))$. \ Since $T$ is
$a$-isoloid, we notice that Weyl's theorem holds for $f(T)$ by the proof of
Theorem \ref{thm25}. \ Therefore $f(T)-\lambda$ has closed range for all
$\lambda\in\pi_{00}^{a}(f(T))$, and it follows from Corollary \ref{cor33} that
$a$-Weyl's theorem holds for $f(T)$.
\end{proof}

\end{document}